\documentclass[letterpaper, 10 pt, conference]{ieeeconf}  % Comment this line out
\IEEEoverridecommandlockouts                              % This command is only
\usepackage{graphics} % for pdf, bitmapped graphics files
\usepackage{epsfig} % for postscript graphics files
\usepackage{amsmath} % assumes amsmath package installed
\usepackage{amssymb}  % assumes amsmath package installed
\usepackage{flushend}

\usepackage{ifthen}
\usepackage{soul}
\newboolean{arxiv}

\usepackage{amsthm}  %Need to check 
\usepackage{mathtools}
\usepackage{commath}
\usepackage{amsfonts}
\usepackage{subcaption}
\usepackage[linesnumbered,ruled,vlined]{algorithm2e} %algorithm package
\SetKwInput{KwInput}{Input}                % Set the Input
\SetKwInput{KwOutput}{Output}
\SetKwInput{KwInitialize}{Initialize} 
\usepackage{xcolor} % I ADDED THIS! Only for proof-reading
\usepackage{ifthen}
\usepackage{soul}

\newtheorem{theorem}{Theorem}

\newtheorem{definition}[theorem]{Definition}

\usepackage[noadjust]{cite}

\title{\LARGE \bf
A Market Mechanism for a Two-stage Settlement Electricity Market with Energy Storage 
}

\author{Rajni Kant Bansal, Enrique Mallada, and Patricia Hidalgo-Gonzalez% <-this % stops a space
\thanks{Rajni Kant Bansal and Patricia Hidalgo-Gonzalez are with the %Center for Energy Research and Department of Mechanical and Aerospace Engineering at the 
University of California San Diego, USA. 
        {Emails: \tt\small \{rabansal, phidalgogonzalez\}@ucsd.edu}}%
\thanks{Enrique Mallada is with the %Center for Energy Research and Department of Mechanical and Aerospace Engineering at the 
Johns Hopkins University, USA. 
        {Email: \tt\small \{mallada\}@jhu.edu}}
}

\begin{document}

\setboolean{arxiv}{true}  

\maketitle
\thispagestyle{empty}
\pagestyle{empty}

%%%%%%%%%%%%%%%%%%%%%%%%%%%%%%%%%%%%%%%%%%%%%%%%%%%%%%%%%%%%%%%%%%%%%%%%%%%%%%%%
\begin{abstract}

% v1
Electricity markets typically clear in two stages: a day-ahead market and a real-time market. In this paper, we propose market mechanisms for a two-stage multi-interval electricity market with energy storage, generators, and demand uncertainties. We consider two possible mixed bidding strategies: storage first bids cycle depths in the day ahead and then charge-discharge power bids in the real-time market for any last-minute adjustments. While the first strategy only considers day-ahead decisions from an individual participant's perspective as part of their individual optimization formulation, the second strategy accounts for both the market operator's and participants' perspectives. We demonstrate that the competitive equilibrium exists uniquely for both mechanisms. However, accounting for the day-ahead decisions in the bidding function has several advantages. Numerical experiments using New York ISO data provide bounds on the proposed market mechanism.
\end{abstract}

%%%%%%%%%%%%%%%%%%%%%%%%%%%%%%%%%%%%%%%%%%%%%%%%%%%%%%%%%%%%%%%%%%%%%%%%%%%%%%%%
\section{Introduction}\label{sec_1}

% v1
Energy storage systems, like grid-scale lithium-ion batteries, that are characterized by a fast response time, high ramp rate, their ability to complement variability in renewable energy dispatch, etc., are being considered across the grid for essential services. Recent policy steps from regulators allowing market participation of emerging technologies, including energy storage, have further increased its adoption in electricity markets~\cite{aeo2021,ferc2222,denholm2019potential}. In particular, several works have investigated the benefits of energy storage for system reliability and power quality, e.g., supporting renewable energy resources, transmission and distribution networks, etc.~\cite{solar_integr,wind_integ,trans_distr_supp,freq_reserve_req,bansal_epri}. However, determining energy storage’s marginal operation costs remains a key challenge in resource participation and efficient market clearing. Unlike conventional fuel-based generators, where marginal costs depend on energy supply, or variable renewable energy resources, with nearly zero operating costs, the degradation incurred due to charge-discharge cycles constitutes the bulk of the operation cost for energy storage. 

Recent efforts have developed participation bids while accounting for the operation cost of storage in two ways. The first approach relies on creating an optimal sequence of charge-discharge energy bids, i.e., price-quantity pairs, typically using a Bilevel optimization framework with the market operator as a follower~\cite{he2015optimal,dheepak_model_1,ramteen_model_1}. However, the resulting market mechanisms based on energy bids require stringent conditions to align with the social optimum and result in incentive-misaligned market outcomes~\cite{bansal2021market}. The second approach relies on developing a state of charge (SoC) based charge-discharge bids to account for SoC-dependent energy storage physical characteristics, e.g., storage degradation~\cite{BXU_model_2, Cong_model_2}. These strategies either assume or estimate unknown prices in markets. Further, the formulations typically consider a single stage, e.g., real-time market, and do not take into account the impact of the decisions in the day-ahead market on the subsequent decisions in the real-time market. Besides, these works provide little insight into the impact of storage decisions on the resulting market equilibrium.

%for the impact of resource owner or market operator decisions in the day-ahead market on the subsequent decisions in the real-time market.   

In this paper, we study energy storage participation in a two-stage multi-interval market with generators under demand uncertainties. We propose mixed bidding market mechanisms with storage bidding an energy-cycling function, i.e., charge-discharge cycle depths as a function of per-cycle prices, while a generator bids a supply function in the day-ahead market. The day-ahead market clears efficiently based on the forecast for the next day. The formulation is implemented via a convex optimization problem that utilizes the convex storage degradation cost function as a combination of the Rainflow cycle counting algorithm with a cycle stress function~\cite{bansal2020storage,bansal2021market}. 

The cost of degradation depends on temporally coupled charge-discharge cycles. Such temporal coupling makes it challenging to design bids for the real-time stage that clears at a relatively faster timescale with adjustments to day-ahead commitments. To address this, we propose two possible participation strategies: 1) Participants account for their day-ahead decisions by including them in their individual problems, i.e., profit maximization. More specifically, we propose a constrained energy bid for storage in the real-time market, such that the storage physical characteristics, e.g., storage degradation function, remain the same. 2) Participants account for their day-ahead market decisions in both the bidding function and individual problems. A generator bids a supply function in the real-time market, similar to the day-ahead stage. The real-time market clearing incorporates a rolling time horizon window that includes binding and advisory intervals. The market operator optimizes the dispatch decision with an updated forecast for the demand. 

\subsubsection*{Contributions} Our work has multiple contributions. \textit{First}, we extend the previously proposed cycle-based market mechanism~\cite{bansal2021market} for the day-ahead market and provide a uniform price market mechanism for storage. The extended design allows market operators to generate a uniform per-cycle price for various storage units, similar to the case of generators in the existing market design. This modification is necessary for the real-time market design. \textit{Second}, we study two mixed market mechanisms and the resulting competitive equilibrium in the real-time market. It assumes participants as price-takers, and we show that the equilibrium exists uniquely. \textit{Third,} we provide a case study for the proposed market mechanism and do a comparative analysis with the social planner problem as we change the physical characteristics of energy storage. Our analysis shows that the proposed mechanism performs within the bounds of the benchmark social planner problem, meaning a bounded loss in the incentive alignment with participation in the real-time market.

The rest of the paper is organized as follows. In section~\ref{sec_2} we introduce the market mechanisms for day-ahead and real-time market. The uniform price market mechanism for the day-ahead market is discussed in Section~\ref{sec_3}. In Section~\ref{sec_4} we compare the two mixed market mechanisms and analyze the existence of competitive equilibrium for the proposed. Section~\ref{sec_5} discusses the case study and Section~\ref{sec_6} concludes the paper.

\paragraph*{Notation} We use $||\cdot||_2^2$ to denote the Euclidean norm and $\langle \cdot, \cdot \rangle$ to denote the inner product. Also, $f(a;b)$ denotes a
function of independent variable $a$ with $b$ as a parameter.

\section{Market Model}\label{sec_2}

In this section, we formulate the social planner problem and describe the proposed market mechanism for a two-stage settlement multi-interval electricity market.

\subsection{Model Preliminaries}
 Consider a two-stage market consisting of a day-ahead and a real-time market, where a set $\mathcal{G}$ of generators and a set $\mathcal{S}$ of storage units participate to meet inelastic demand over a time horizon $t \in \mathcal{T}:= \{1,..., T\}$. We denote the demand forecast in day-ahead as $d^d \in \mathbb{R}^T$ and the forecast error in real-time as $d^r \in \mathbb{R}^T$. Thus, the net demand $d \in \mathbb{R}^T$ over two stages is given by
\begin{align}\label{demand_notation}
    d:= d^d+d^r
\end{align}
Each generator $j \in \mathcal{G}$ dispatches $g_j^d \in \mathbb{R}^T$ and $g_j^r\in \mathbb{R}^T$ over the time horizon $\mathcal{T}$ in the day-ahead and real-time markets, respectively. For $v \in \{d,r\}$, the stage-wise and overall power dispatch of generator $j$ is subject to capacity constraints 
\begin{subequations}\label{gen_ineq_constraint}
\begin{align}
    &\underline{g}_j \le g_{j,t}^v \le \overline{g}_j,  \ t \in \mathcal{T} \label{gen_cap_individual} \\
    & \underline{g}_j \le g_{j,t}^d +  g_{j,t}^r \le \overline{g}_j,  \ t \in \mathcal{T} \label{gen_cap_total}
\end{align}
\end{subequations}
where $\underline{g}_j,\overline{g}_j$ denote the minimum and maximum generation limits, respectively. Similarly, each energy storage $s \in \mathcal{S}$ of capacity $E_s$ dispatches $u_s^d \in \mathbb{R}^T \text{and } u_s^r\in \mathbb{R}^T$ in the day-ahead and real-time markets, respectively. The individual and total energy dispatch, i.e., discharge (positive) or charge (negative) rate, is bounded as,
\begin{subequations}\label{str_ineq_constraint}
\begin{align}
    &\underline{u}_s \le u_{s,t}^v \le \overline{u}_s,  \ t \in \mathcal{T} \label{str_cap_individual} \\
    & \underline{u}_s \le u_{s,t}^d +  u_{s,t}^r \le \overline{u}_s,  \ t \in \mathcal{T} \label{str_cap_total}
\end{align}
\end{subequations}
where $v \in \{d,r\}$. The pair $\underline{u}_s,\overline{u}_s$ denote the minimum and maximum storage rate limits, respectively. The net output of generator $j$ and storage $s$ over two stages is denoted as,
\begin{align} \label{net_dispatch_notation}
    g_j := g_j^d+g_j^r, j\in \mathcal{G}, \ u_s:= u_s^d + u_s^r, s \in \mathcal{S}
\end{align}
We denote the stored energy in storage $s$ by a normalized State of Charge (SoC) profile $x_s^d \in \mathbb{R}^{T+1}$ with initial SoC $x_{s,0}^d$ associated with the day-ahead storage dispatch $u_s^d$. The SoC evolves over time horizon $\mathcal{T}$ as, 
\begin{align}\label{SoC_evolve_da}
        x_{s,t}^d - x_{s,t-1}^d = \frac{1}{E_s}u_{s,t}^d, \forall s \in \mathcal{S}, t \in \mathcal{T}
\end{align}
Analogously, we denote SoC profile $x_s\in \mathbb{R}^{T+1}$ with initial SoC $x_{s,0}$ corresponding to the net storage dispatch $u_s$, such that the SoC evolves as
\begin{align}\label{SoC_evolve_all}
        x_{s,t} - x_{s,t-1} = \frac{1}{E_s}u_{s,t}, \forall s \in \mathcal{S}, t \in \mathcal{T}
\end{align}

The normalized SoC is bounded as,
\begin{subequations}\label{SoC_ineq_all}
\begin{align}
     & 0 \le x_{s,t}^d \le 1, s \in \mathcal{S}, t\in \mathcal{T} \label{SoC_ineq_da}\\
     & 0 \le x_{s,t} \le 1, s \in \mathcal{S}, t\in \mathcal{T} \label{SoC_ineq_total}
\end{align}
\end{subequations}

Lastly, we assume periodicity constraints on storage to account for its cyclic nature in current markets\footnote{However, we relax such a periodicity constraint in real-time to allow storage for any immediate adjustments for system reliability.}, i.e., 
\begin{subequations}\label{SoC_periodicity_temp_all}
\begin{align}%\label{SoC_periodicity_temp_da}
        & x_{s,0}^d = x_{s,T}^d \\
        & x_{s,0} = x_{s,T}
\end{align}
\end{subequations}
Substituting~\eqref{SoC_evolve_all} in \eqref{SoC_periodicity_temp_all}, we get
\begin{subequations}\label{SoC_periodicity_all}
\begin{align}
        & \sum_{t\in \mathcal{T}} u_{s,t}^d = 0, \ s \in \mathcal{S} \label{SoC_periodicity_da}\\
        & \sum_{t\in \mathcal{T}} u_{s,t} = 0, \ s \in \mathcal{S} \label{SoC_periodicity_total}
\end{align}
\end{subequations}

\subsection{Social Planner}
The social planner's problem, assuming perfect foresight is given by
\begin{subequations}\label{social_planner}
\begin{align}
    \min_{g_j, j\in \mathcal{G}, u_s, s \in \mathcal{S}, x_s, s\in \mathcal{S}} & \ \sum_{j \in \mathcal{G}} C_j(g_j) + \sum_{s \in \mathcal{G}} C_s(u_s)\\
    & \textrm{s.t.} \sum_{j\in\mathcal{G}}g_j+\sum_{s\in\mathcal{S}}u_s = d \label{two_stage_power_balance}\\
    & \eqref{gen_cap_total},\eqref{str_cap_total},\eqref{SoC_evolve_all},\eqref{SoC_ineq_total},\eqref{SoC_periodicity_total} \nonumber
\end{align}
\end{subequations}
where~\eqref{two_stage_power_balance} denotes the power balance constraint over two stages. In this paper, we assume a quadratic cost function for the generators, given by
\begin{align} \label{gen_cost_func}
    C_j(g_j) = \frac{c_j}{2}\sum_{t\in\mathcal{T}}g_{j,t}^2 + a_j\sum_{t\in\mathcal{T}}g_{j,t}
\end{align}
where $c_j,a_j$ are the cost coefficients\footnote{For ease of analysis, we assume that $a_j =0$. However, the analysis is generalizable for the case $a_j\neq 0$.}. For storage $s$, we adopt the convex cycle-based degradation cost as its operational cost. It combines the Rainflow cycle counting algorithm with a cycle stress function to identify and penalize the cost of charge-discharge cycles~\cite{bansal2021market}. The cost function is given by
\begin{align}\label{str_cost_func}
    C_s(u_s) = \frac{b_s}{2}\nu_s^T\nu_s = \frac{b_s}{2}u_s^TN(u_s)^TN(u_s)u_s 
\end{align}
where $b_s$ is the empirical cost coefficient. Also, $\nu_s \in \mathbb{R}^T$ represents the vector of half-cycle depths created from the Rainflow algorithm, which maps a storage dispatch profile to associated charge-discharge half-cycles, such that 
\begin{align}\label{nu_relation}
    \nu_s := Rainflow(u_s) = N(u_s)u_s
\end{align}
Here the matrix $N(u_s)$ is a function of the storage dispatch $u_s$ and represents the Rainflow map between the storage dispatch profile and half-cycle depths. See, e.g.,~\cite{bansal2021market} for more details.

%PHG continue here
\subsection{Two-Stage Market Mechanism}

\begin{figure}[tbp]
    \centering
    {\includegraphics[width=0.9\linewidth]{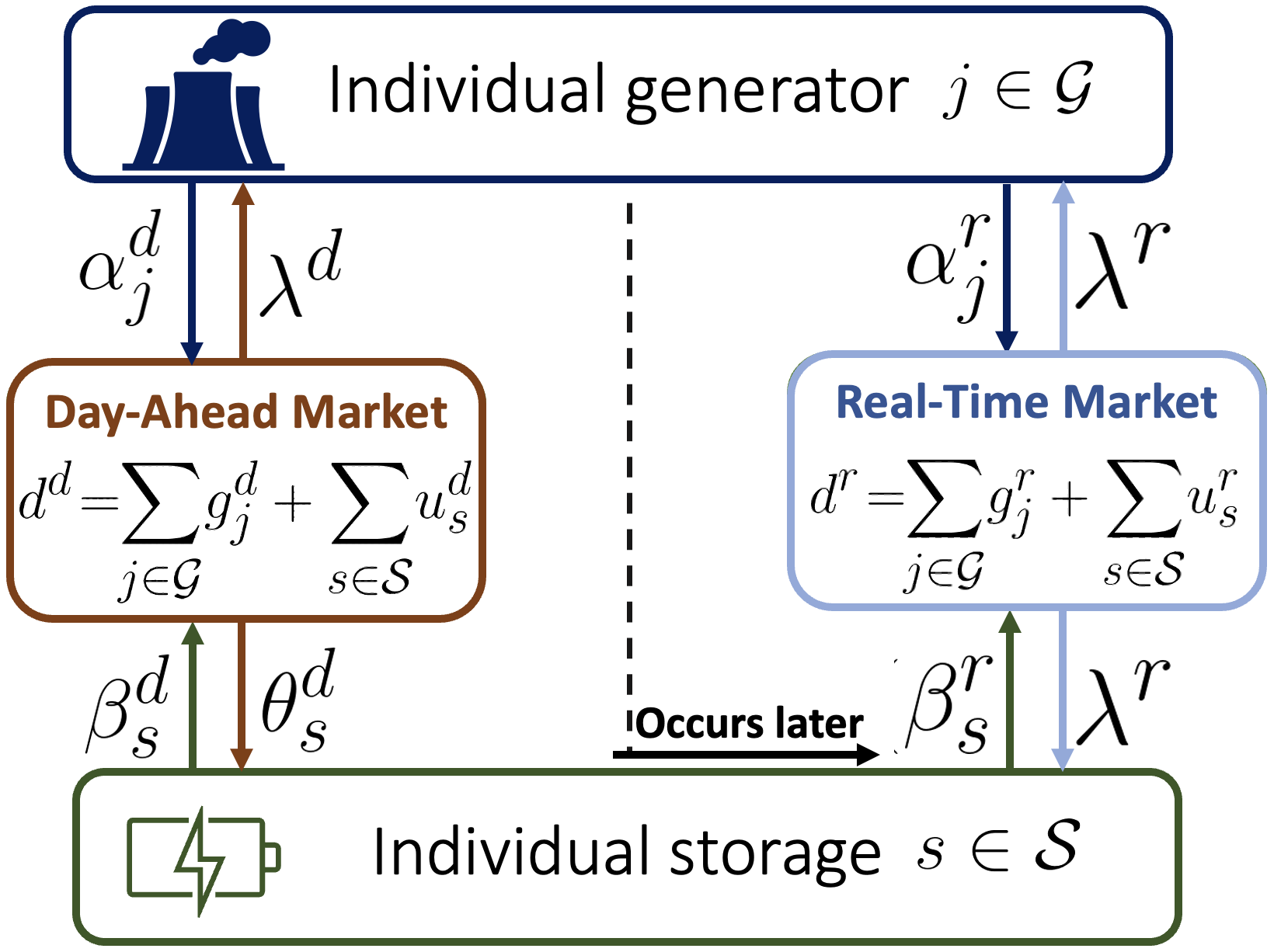} }
    \caption{Two-stage Mixed Market Mechanism.}%
    \label{fig:two_stage}%
\end{figure}

In this subsection, we describe the two-stage market clearing, as shown in Figure~\ref{fig:two_stage}. 
\subsubsection{Day-Ahead Market}: Each generator $j$ submits a supply function parameterized by $\alpha_j^d\in\mathbb{R}$ indicating its willingness to participate in day-ahead market as,
\begin{align}\label{gen_bid_da}
    g_j^d = \alpha_j^d\lambda^d
\end{align}
where $\lambda^d \in \mathbb{R}^T$ denotes the clearing prices in the day ahead. Analogously, each storage $s$ submits an energy cycling function parameterized by $\beta_s^d \in \mathbb{R}$ that maps cycle depths $\nu_s^d \in \mathbb{R}^T$ to per cycle prices $\theta_s^d \in \mathbb{R}^T$, as
\begin{align}\label{str_bid_da}
    \nu_s^d = \beta_s^d\theta_s^d
\end{align} 
The market operator collects all the bids and solves the day-ahead market clearing problem, given by 
\begin{subequations}\label{da_mkt_clearing}
\begin{align}
    \min_{g_j^d, j\in \mathcal{G},(u_s^d,\nu_s^d), s\in\mathcal{S}} & \ \sum_{j\in \mathcal{G}}\sum_{t\in \mathcal{T}} \frac{(g_{j,t}^d)^2}{2\alpha_j} \!+\!\! \sum_{s\in \mathcal{S}}\sum_{t\in \mathcal{T}} \frac{(\nu_{s,t}^d)^2}{2\beta_s}\!\\
  \textrm{s.t.}   &  \sum_{j\in\mathcal{G}}g_j^d + \sum_{s\in\mathcal{S}}u_s^d = d^d \label{da_power_balance}\\
    & \nu_s^d = N(u_s^d)u_s^d \label{rainflow_const_da}\\
    & \eqref{gen_cap_individual},\eqref{str_cap_individual},\eqref{SoC_periodicity_da} \nonumber
\end{align}
\end{subequations}
where~\eqref{rainflow_const_da} denotes the constraint associated with the Rainflow algorithm. The market clearing gives the optimal dispatch and prices, such that generator $j$ produces $g_j$ and gets paid $(\lambda^d)^Tg_j^d$ while storage $s$ produces a cycle depth schedule $\nu_s$ and gets paid $(\theta_s^d)^T\nu_s^d$. Here the prices $\lambda^d, \theta_s^d$ are given by the dual variables corresponding to the constraints~\eqref{da_power_balance} and \eqref{rainflow_const_da}, respectively.

In the day-ahead market, participants are competing against each other to maximize their profit. We assume participants as price-takers such that the individual problems of generator $j$, is
\begin{align}\label{da_gen_prft}
    \max_{g_j^d} \ \langle\lambda^d, g_j^d\rangle - C_j(g_j^d), \ \textrm{s.t.} \ \eqref{gen_bid_da}
\end{align}
% , g_{j,t} \in (\underline{g}_j,\overline{g}_j), t\in \mathcal{T}
and storage $s$, is
\begin{align}\label{da_str_prft}
    \max_{(u_s^d,\nu_s^d) } \langle\theta_s^d,\nu_s^d\rangle - C_s(u_s^d) \ \textrm{s.t.} \ \eqref{str_bid_da}
\end{align}

Since the participants are price-takers, they do not anticipate their decisions in the real-time market.

\subsubsection{Real-Time Market}
Analogous to the day-ahead market, each generator $j$ submits a supply function $f:\mathbb{R}\times \mathbb{R}^T \rightarrow \mathbb{R}^T$,  parameterized with $\alpha_j^r \in \mathbb{R}$, as
\begin{align}\label{gen_bid_rt_func}
    g_j^r = f(\alpha_j^r;\lambda^r)
\end{align}
where $\lambda^r \in \mathbb{R}^T$ denotes the clearing prices in the real-time. However, unlike storage in day-ahead markets, in the real-time market storage $s$ bids a supply function  $h:\mathbb{R}\times \mathbb{R}^T \rightarrow \mathbb{R}^T$, parameterized by $\beta_s^r$, as 
\begin{align}\label{str_bid_rt_func}
    u_s^r = h(\beta_s^r;\lambda^r)
\end{align}
Given the bids $(\alpha_j^r,\beta_s^r)$ the market operator clears the market and meets the supply-demand balance for each time period $t \in \mathcal{T}$, as
\begin{align}\label{rt_mkt_clearing}
    d^r = \sum_{j\in\mathcal{G}}f(\alpha_j^r;\lambda^r) + \sum_{s\in\mathcal{S}}h(\beta_s^r;\lambda^r)
\end{align}
Given the prices in the real-time market and the participant's decision in the day-ahead market, the individual problem of generator $j$ is given by,
\begin{align}\label{rt_gen_prft}
    \max_{g_j^r} \ \langle\lambda^d, g_j^d \rangle+\langle\lambda^r, g_j^r\rangle - C_j(g_j^d+g_j^r) \ \textrm{s.t.} \ \eqref{gen_bid_rt_func}
\end{align}
and storage $s$ is 
\begin{align}\label{rt_str_prft}
    \max_{u_s^r} \ \langle \theta_s^d , \nu_s^d \rangle +  \langle\lambda^r, u_s^r \rangle - C_s(u_s^d+u_s^r)\ \textrm{s.t.} \ \eqref{str_bid_rt_func}
\end{align}

\subsection{Market Equilibrium}

In this subsection, we characterize the properties of competitive equilibrium in a two-stage market, such that the market clears and participants do not deviate from their bid. 
\begin{definition}
    Each stage of a two-stage market is at the competitive equilibrium if the participant bids and the clearing prices, i.e., $(\lambda^d,\alpha_j^d, \beta_s^d)$ and $(\lambda^r,\alpha_j^r, \beta_s^r)$, in day-ahead and real-time markets, respectively, satisfy:
    \begin{enumerate}
        \item The bid $\alpha_j^d(\alpha_j^r)$ of generator $j$ in the day-ahead (real-time) market maximizes its profit.
        \item The bid $\beta_s^d(\beta_a^r)$ of storage $s$ in the day-ahead (real-time) market maximizes its profit.
        \item The inelastic demand $d^d, d^r$ in the day-ahead and the real-time market is met, resulting in clearing prices $\lambda^d, \theta_s^d$ and $\lambda^r$, respectively. 
    \end{enumerate}
\end{definition}

\section{Uniform Price Market Mechanism}\label{sec_3}

In this section, we extend the previous cycle-aware market mechanism~\cite{bansal2021market} in the day-ahead market and propose a uniform price market mechanism.

\begin{theorem}(Incentive Compatibility~\cite{bansal2021market})
    The competitive equilibrium in the day-ahead market also aligns with the social optimum. 
\end{theorem}

The incentive compatibility results first appeared in~\cite{bansal2021market} for a market mechanism where storage bids cycle depths and generator bids power in a day-ahead market. However, the resulting competitive equilibrium provides non-uniform prices for storage units, which may not be desirable from a market operator perspective. To address this, we propose the following theorem that ascertains uniform prices in a market.

%\subsection{Uniform Price Market Mechanism}

%In this paper, we extend the previous work~\cite{bansal2021market} to a uniform price market mechanism for energy storage in the day-ahead market.

\begin{theorem}
    Let's assume that the capacity constraints of storage~\eqref{str_cap_individual} and generator~\eqref{gen_cap_individual} are non-binding. Furthermore, the SoC constraint associated with energy storage~\eqref{SoC_ineq_da} is also non-binding. Given optimal bids $\alpha_j^d, \beta_s^d$, there exists a unique set of coefficients $\epsilon_s$ such that the day-ahead market clearing~\eqref{da_mkt_clearing} results in uniform prices. %and each generator and storage is paid uniformly with uniform prices $\lambda^d, \theta^d$, respectively.
\end{theorem}

\ifthenelse{\boolean{arxiv}}{\begin{proof}
    Let's assume a proportional energy storage dispatch of the form
\begin{align}\label{str_prop_disp}
    u_s^d = \epsilon_s (d^d-\sum_{j\in\mathcal{G}}g_j^d), \textrm{ where } \sum_{s\in\mathcal{S}} \epsilon_s = 1
\end{align}
Let's denote $u^d:= (d^d-\sum_{j\in\mathcal{G}}g_j^d)$. The KKT conditions of the convex day-ahead market clearing problem~\eqref{da_mkt_clearing} are given by, 
\begin{subequations}\label{kkt_cond_da_mkt}
\begin{align}
    & d^d = \sum_{j\in\mathcal{G}}g_j^d + \sum_{s\in\mathcal{S}}u_s^d , \ g_j^d = \alpha_j^d \lambda^d \label{kkt_cond_da_mkt.a}\\ 
    & \nu_s^d = N(u_s^d)u_s^d, \ \nu_s^d = \beta_s^d\theta_s^d \label{kkt_cond_da_mkt.b}\\
    & \lambda^d = \sum_k\gamma_{k}N_{k}(u^d_s)^TN_k(u_s^d)u_s^d + \delta_s\mathbf{1} \label{kkt_cond_da_mkt.c}\\
    & \mathbf{1}^Tu_s^d = 0. \label{kkt_cond_da_mkt.d}
\end{align}
\end{subequations}
where $\gamma_k$ are convex coefficients associated with subgradients of the piecewise linear convex cost of storage cycling~\eqref{str_cost_func}. See e.g.,~\cite{bansal2021market}, for more details. Also, $\lambda^d, \theta_s^d,$ and $\delta_s$ are the dual variables associated with constraints~\eqref{da_power_balance}, \eqref{rainflow_const_da}, and \eqref{SoC_periodicity_da}, respectively. Substituting the proportional dispatch condition~\eqref{str_prop_disp} in the KKT conditions~\eqref{kkt_cond_da_mkt}, we get
\begin{subequations}\label{kkt_cond_da_mkt_prop}
\begin{align}
    & u^d = d^d - \sum_{j\in\mathcal{G}}g_j^d, \ g_j^d = \alpha_j^d \lambda^d \\ 
    & \nu_s^d = \epsilon_sN(u^d)u^d \\
    & \nu_s^d = \beta_s^d\theta_s^d \implies \theta_s^d = \frac{\epsilon_s}{\beta_s}N(u^d)u^d \\
    & \lambda^d = \sum_k\gamma_{k}\frac{\epsilon_s}{\beta_s}N_{k}(u^d)^TN_k(u^d)u^d + \delta_s\mathbf{1}\\
    & \mathbf{1}^Tu^d = 0
\end{align}
\end{subequations}
 where we use the~\cite[Lemma 1]{bansal2021market} to write 
 \[
    N(u_s^d) = N(\epsilon_su^d) = N(u^d).
 \] 
 Rewriting the KKT conditions~\eqref{kkt_cond_da_mkt_prop}, we get 
\begin{subequations}\label{kkt_cond_da_mkt_prop_rewrite}
\begin{align}
    & \frac{d-u^d}{\sum\limits_{j\in\mathcal{G}}\alpha_j^d} = \sum_k\gamma_{k}\frac{\epsilon_s}{\beta_s^d}N_{k}(u^d)^TN_k(u^d)u^d + \delta_s\mathbf{1} \\
    & \mathbf{1}^Tu^d = 0
\end{align}
\end{subequations}
Given the optimal bids ($\alpha_j^d, \beta_s^d$), we consider a convex optimization problem, as
\begin{subequations}\label{unif_con_obj}
\begin{align}
    \min_{u^d} & \ \frac{|| u^d ||_2^2}{2\sum\limits_{j\in\mathcal{G}}\alpha_j^d} \!-\! \frac{\ \langle d^d,u^d\rangle}{\sum\limits_{j \in \mathcal{G}}\alpha_j^d} \!+\! \frac{|| N(u^d)u^d||_2^2}{\sum\limits_{s \in \mathcal{S}}\beta_s^d} \\
    \textrm{s.t.} & \ \mathbf{1}^Tu^d = 0 :\delta
\end{align}
\end{subequations}
From~\cite{bansal2020storage}, we observe that the optimization problem~\eqref{unif_con_obj} is convex. Therefore, $\exists $ a unique primal-dual solution $(u^d)^{*},\delta^{*}$ to the convex optimization problem~\eqref{unif_con_obj} which satisfies the KKT conditions ~\eqref{kkt_cond_da_mkt_prop_rewrite}. Hence, assuming 
\begin{align}
    \epsilon_s = \frac{\beta_s}{\sum\limits_{s\in\mathcal{S}}\beta_s^d}, \ \delta_s^{*}:= \delta^{*} \ \forall s \in \mathcal{S}, \ u_s^d = \epsilon_su^d
\end{align}
we observe that we have a unique optimal solution with uniform clearing prices 
\begin{subequations}
\begin{align}
    & \lambda^d = \frac{d-u^d}{\sum\limits_{j\in\mathcal{G}}\alpha_j^d}, \\
    & \theta_s^d = \theta^d := \frac{\epsilon_s}{\beta_s^d}N(u^d)u^d = \frac{1}{\sum\limits_{s\in \mathcal{S}}\beta_s^d}N(u^d)u^d
\end{align}
\end{subequations}
that satisfies the KKT conditions of the original problem. Therefore, we have a uniform price market mechanism.
\end{proof}
}
{{The proof is provided in \addcite.}} 

For ease of analysis, we assume the inequality constraints are non-binding. However, the result can be observed for a general formulation as well. Moreover, the resulting competitive equilibrium for the uniform price market mechanism also aligns with the social planner optimum in the day-ahead market.

\section{Competitive Equilibrium in the Real-Time Market}\label{sec_4}

In this section, we discuss two potential energy storage participation strategies in the real-time market, i.e., day-ahead unaware bidding and day-ahead aware bidding. For ease of exposition and closed-form analysis, we consider in this section a simplified setting where the power balance~\eqref{rt_mkt_clearing} is the only constraint in the real-time market clearing. A general setting is considered in the case study in Section~\ref{sec_5}. 

\subsection{Mixed Market Mechanism - Day-Ahead Unaware Bidding} \label{subsec_2}

In this subsection, we characterize the competitive equilibrium in the real-time market given the participants' decision in the day-ahead market. Each generator $j$ bids 
\begin{align}\label{gen_bid_rt}
    g_j^r = \alpha_j^r\lambda^r
\end{align}
and storage $s$ bids 
\begin{align}\label{str_bid_rt}
    u_s^r = \beta_s^r\lambda^r
\end{align}
in the real-time market. Substituting \eqref{gen_bid_rt} in \eqref{rt_gen_prft}, we get
\begin{align}\label{rt_gen_prft_bid}
    \max_{\alpha_j^r} \ \alpha_j^r||\lambda^r||_2^2 - \frac{c_j}{2}||g_j^d+\alpha_j^r\lambda^r||_2^2
\end{align}
where we drop the revenue term $\langle\lambda^d,g_j^d\rangle$ from the day-ahead market as it does not affect the objective in the optimization problem. Similarly, substituting~\eqref{str_bid_rt} in \eqref{rt_str_prft}, we get
\begin{align}\label{rt_str_prft_bid}
    \max_{\beta_s^r} \ \beta_s^r||\lambda^r||_2^2 - \frac{b_s}{2}||N(u_s^d+\beta_s^r\lambda^r)(u_s^d+\beta_s^r\lambda^r)||_2^2
\end{align}

%For ease of exposition and closed-form analysis, we consider in this section a simplified setting where the power balance~\eqref{rt_mkt_clearing} is the only constraint in the real-time market clearing. A general setting is considered in the case study in the Section~\ref{sec_4}. 

To account for the temporally coupled cost of storage degradation while making decisions in the real-time market at a faster time scale, we assume that storage deviations are constrained in the neighborhood of its day-ahead decision such that
\begin{align}
    N(u_s^d+u_s^r) = N(u_s^d)
\end{align}
More precisely, we constrain the storage dispatch such that the piecewise linear map between cycle depths and SoC profile remains the same. Such a strategy allows storage to participate in the market with a simplified setting, e.g., existing supply function, with the flexibility to change the optimal bid to account for the degradation cost. %We also propose an algorithm in the Appendix that will enable us to constrain the energy storage in real-time market.

\begin{theorem}
    Let's assume a uniform price market mechanism in the day-ahead market. Also, assume that 
    \[
        N(u_s^d+u_s^r) = N(u_s^d).
    \]
    The competitive equilibrium in the real-time market, given the participants' decision in the day-ahead market, exists uniquely. There exists a unique set of coefficients $\omega$ such that the equilibrium can be determined as
    \begin{subequations}\label{rt_comp_eqbm}
    \begin{align}
        & \lambda_r = \omega d^r\\
        & \alpha_j^r = c_j^{-1} - \frac{<g_j^d,\lambda^r>}{||\lambda^r||_2^2} \\
        & \beta_s^r = b_s^{-1}\frac{||\lambda^r||_2^2 - \langle \theta_s^d,N(u_s^d)\lambda^r\rangle}{||N(u_s^d)\lambda^r||^2_2} \\
        & \!\!\omega \!= \!\frac{<\!\lambda^d,d^r\!\!>}{||d^r||_2^2} \!+\!\!\! {\left(\!\sum\limits_{j\in\mathcal{G}}\!\!c_j^{-1}\!\!+\!\!\!\sum\limits_{s \in \mathcal{S}}\!b_s^{-1}\frac{||d^r||_2^2}{||\!\sum\limits_k\!\!\gamma_kN(\!u_s^d)d^r||^2_2}\!\!\right)\!\!\!}^{-1}
    \end{align}
    \end{subequations}
    where $\gamma_k \ge 0, \sum_k \gamma_k = 1$ are the convex coefficients associated with the subgradients of the piecewise linear convex cost function~\eqref{str_cost_func}. See e.g., \cite{bansal2021market} for more details.  
\end{theorem}

\ifthenelse{\boolean{arxiv}}{\begin{proof}
    Given the day-ahead decisions and prices in the real-time market, we can solve for the optimal bid by taking the derivative of the convex optimization problem~\eqref{rt_gen_prft_bid}, as
    \begin{subequations}\label{kkt_gen_bid_rt}
    \begin{align}
    \!\!\!\!||\lambda^r||^2_2 -\! c_j(g_j^d\!+\!\alpha_j^r\lambda^r)^T\!\lambda^r \!\!=\! 0 \!\!
    \implies \!\! \alpha_j^r \!=\! c_j^{-1} \!-\! \frac{\langle g_j^d, \lambda^r\rangle}{||\lambda^r||^2_2} \label{kkt_gen_bid_rt.b}
    \end{align}
    \end{subequations}
    Similarly, taking the derivative of \eqref{rt_str_prft_bid}, we get the optimal bid of energy storage $s$,
    \begin{align}\label{kkt_str_bid_rt_tmp}
    & \!||\lambda^r||^2_2 - b_s\!\left(N(u_s^d)(u_s^d+\beta_S^r\lambda^r)\!\right)^T\!N(u_s^d)\lambda^r \!=\! 0
    \end{align}
    where we use $N(u_s^d+u_s^r) = N(u_s^d)$. Substituting~\eqref{nu_relation} and \eqref{str_bid_da} in \eqref{kkt_str_bid_rt_tmp}, we get
    \begin{subequations} \label{kkt_str_bid_rt}
    \begin{align}
    \implies & ||\lambda^r||^2_2 - b_s\left(\nu_s^d+\beta_S^rN(u_s^d)\lambda^r\right)^TN(u_s^d)\lambda^r \!=\! 0 \label{kkt_str_bid_rt.a}\\
    \implies & ||\lambda^r||^2_2 \!-\! \left(\theta_s^d+b_s\beta_S^rN(u_s^d)\lambda^r\right)^T\!\!N(u_s^d)\lambda^r \!=\! 0 \label{kkt_str_bid_rt.b} \\
    \implies & \beta_S^r = b_s^{-1}\frac{||\lambda^r||_2^2 - \langle \theta_s^d,N(u_s^d)\lambda^r\rangle}{||N(u_s^d)\lambda^r||^2_2} \label{kkt_str_bid_rt.c}
    \end{align}
    \end{subequations}
    At the equilibrium, \eqref{rt_mkt_clearing}, \eqref{kkt_gen_bid_rt.b}, and \eqref{kkt_str_bid_rt.c} must hold simultaneously. Since $\lambda^r$ is proportional to $d^r$, let's assume that $\exists \ \omega\in\mathbb{R}$ such that $\lambda^r = \omega d^r$. Substituting~\eqref{gen_bid_da} in \eqref{kkt_gen_bid_rt.b}, as
    \begin{subequations}\label{kkt_gen_bid_rt_v2}
    \begin{align}
        & \!\!\alpha_j^r \!=\! c_j^{-1}\!\left(\!1 \! -\! \frac{<\lambda^d,\lambda^r>}{||\lambda^r||_2^2}\right) \!\!=\! c_j^{-1}\!\left(\!\!1 \!-\! \frac{1}{\omega}\frac{<\lambda^d,d^r>}{||d^r||_2^2}\!\!\right) \label{kkt_gen_bid_rt_v2.a}\\
        & \implies \sum\limits_{j \in \mathcal{G}}\alpha_j^r = \sum\limits_{j \in \mathcal{G}}c_j^{-1}\left(1 - \frac{1}{\omega}\frac{<\lambda^d,d^r>}{||d^r||_2^2}\right)\label{kkt_gen_bid_rt_v2.b}
    \end{align}
    \end{subequations}
Similarly, substituting \eqref{kkt_cond_da_mkt.c} in  \eqref{kkt_str_bid_rt.c} and simplifying as
    \begin{align}
         \!\!\!\!\!\!\beta_s^r \!=\! b_s^{-1}\frac{||\lambda^r||_2^2 - \langle \lambda^d,\lambda^r\rangle}{||\!\sum\limits_k\!\gamma_kN(u_s^d)\lambda^r||^2_2} \!= b_s^{-1}\frac{||d^r||_2^2 - \omega^{-1}\langle \lambda^d,d^r\rangle}{||\sum\limits_k\!\gamma_kN(u_s^d)d^r||^2_2} \label{kkt_str_bid_rt_v2.b} 
    \end{align}
    Substituting~\eqref{kkt_gen_bid_rt_v2.b} in \eqref{kkt_str_bid_rt_v2.b}, as
    \begin{subequations}\label{kkt_str_bid_rt_v3}
    \begin{align}
         \implies & \beta_s^r = b_s^{-1}\frac{\sum\limits_{j\in\mathcal{G}}\alpha_j}{\sum\limits_{j\in\mathcal{G}}c_j^{-1}}\frac{||d^r||_2^2}{||\sum_k\gamma_kN(u_s^d)d^r||^2_2} \label{kkt_str_bid_rt_v3.a}\\
         \implies & \sum\limits_{s\in\mathcal{S}}\beta_s^r = \sum\limits_{s\in\mathcal{S}}b_s^{-1}\frac{\sum\limits_{j\in\mathcal{G}}\alpha_j}{\sum\limits_{j\in\mathcal{G}}c_j^{-1}}\frac{||d^r||_2^2}{||\sum_k\gamma_kN(u_s^d)d^r||^2_2} \label{kkt_str_bid_rt_v3.b}
    \end{align}
    \end{subequations}
    where we use the fact that energy dispatch is proportional and $N(u_s^d) = N(u^d)$ in the case of uniform price market mechanism in equation~\eqref{kkt_str_bid_rt_v3.b}. Substituting \eqref{kkt_gen_bid_rt_v2.b} and \eqref{kkt_str_bid_rt_v3.b} in \eqref{rt_mkt_clearing}, we get
    \begin{align}
        \omega \!= \!\frac{<\!\!\lambda^d\!,d^r\!\!>}{||d^r||_2^2} \!+\!\! \left(\!\sum\limits_{j\in\mathcal{G}}\!\!c_j^{-1}\!\!+\!\!\!\sum\limits_{s \in \mathcal{S}}\!b_s^{-1}\frac{||d^r||_2^2}{||\sum\limits_k\gamma_kN(u_s^d)d^r||^2_2}\!\!\right)^{-1} \nonumber
    \end{align}    
    Hence, under the price-taking assumption the optimal bid and clearing prices $(\alpha_j^r, j \in \mathcal{G}, \beta_s^r, s\in \mathcal{S}, \lambda^r)$ exist uniquely. 
\end{proof}
}
{{The proof is provided in \addcite.}} 

Although the competitive equilibrium exists uniquely, it requires a customized iterative algorithm, where participants' bids and clearing prices are updated repetitively until convergence, to solve the equilibrium, which may not be feasible for real-time market operations. To address this, we propose a modified bidding approach in the following subsection that accounts for the day-ahead decisions in the bidding function. 

\subsection{Mixed Market Mechanism - Day-Ahead Aware Bidding}\label{subsec_3}

Given the day-ahead decisions, each generator $j$ bids 
\begin{align}\label{gen_bid_rt_v2}
    g_j^r = \alpha_j^r\lambda^r - g_j^d
\end{align}
and storage $s$ bids 
\begin{align}\label{str_bid_rt_v2}
    u_s^r = \beta_s^r\lambda^r - u_s^r
\end{align}
in the real-time market. Substituting \eqref{gen_bid_rt_v2} in \eqref{rt_gen_prft}, we get,
\begin{align}\label{rt_gen_prft_bid_v2}
    \max_{\alpha_j^r} \ \alpha_j^r||\lambda^r||_2^2 + \langle \lambda^r, g_j^d \rangle - \frac{c_j}{2}||\alpha_j^r\lambda^r||_2^2
\end{align}
and substituting~\eqref{str_bid_rt_v2} in \eqref{rt_str_prft}, we get
\begin{align}\label{rt_str_prft_bid_v2}
    \max_{\beta_s^r} \ \beta_s^r||\lambda^r||_2^2 + \langle \lambda^r, u_s^d \rangle - \frac{b_s}{2}||N(\beta_s^r\lambda^r)(\beta_s^r\lambda^r)||_2^2
\end{align}
Note that we drop the terms associated with the day-ahead revenue in \eqref{rt_gen_prft_bid_v2} and \eqref{rt_str_prft_bid_v2}, as it does not affect the objective in the individual problem. Once all the bids $(\alpha_j^r, \beta_s^r)$ are collected, the market operator clears the real-time market to achieve supply-demand balance, given by
\begin{align}\label{rt_mkt_clearing_v2}
    d^r = \sum_{j\in\mathcal{G}}\left(\alpha_j^r\lambda^r - g_j^d\right)+ \sum_{s\in\mathcal{S}}\left(\beta_s^r\lambda^r - u_s^d\right)
\end{align}

\begin{theorem}\label{mix_thrm_v2}
    The competitive equilibrium of the mixed market mechanism with day-ahead decision-aware bids in the real-time market exists uniquely, as 
    \begin{subequations}
    \begin{align}
    & \lambda^r = \phi d, \\
    &\alpha_j^r = c_j^{-1},\\
    &{\beta_s^r} = b_s^{-1}\frac{||\lambda^r||_2^2}{||N(\lambda^r)\lambda||_2^2} \\
    &\phi^{-1}  = \left(\sum_{i\in \mathcal{S}}b_s^{-1}\frac{||d||_2^2}{||N(d)d||_2^2}+\sum_{j\in \mathcal{G}}c_j^{-1} \right)
    \end{align}
    \end{subequations}
\end{theorem}

\ifthenelse{\boolean{arxiv}}{\begin{proof}
    Writing the derivative of \eqref{rt_gen_prft_bid_v2} for the optimal bid of each generator $j$, we get
    \begin{align}
    & ||\lambda^r||_2^2 ( 1 - \alpha_j^rc_j)= 0  \implies {\alpha_j^r} = c_j^{-1} \label{opt_alpha}, \ \forall j \in \mathcal{G}
\end{align}

Similarly, we take the derivative of \eqref{rt_str_prft_bid_v2} for storage, as
\begin{subequations}
\begin{align}
 \!\!\!\!||\lambda^r||_2^2 \!-\! b_s{\beta_s^r}||\!N(\!\lambda^r\!)\lambda^r||_2^2 \!=\! 0 \!\!\!
    \implies\!\!\! {\beta_s^r} \!\!=\! b_s^{-1}\frac{||\lambda^r||_2^2}{||\!N(\!\lambda^r\!)\lambda^r||_2^2} \label{opt_beta_hat}
\end{align}
\end{subequations}
where we use ~\cite[Lemma 1]{bansal2021market} in \eqref{opt_beta_hat}. At the equilibrium ~\eqref{rt_mkt_clearing_v2}, \eqref{opt_alpha}, and \eqref{opt_beta_hat} must hold simultaneously. Since $\lambda^r$ is proportional to $d$, let's assume $\exists \ \phi \in \mathbb{R}$ such that $\lambda^r = \phi d$. Substituting \eqref{da_power_balance}, \eqref{opt_alpha}, \eqref{opt_beta_hat} in \eqref{rt_mkt_clearing_v2}, we can solve for $\phi$ as
\begin{align}
     \phi^{-1} =\left(\sum_{s\in \mathcal{S}}b_s^{-1} \frac{||d||_2^2}{||N(d)d||_2^2}+\sum_{j\in \mathcal{G}}{c_j}^{-1} \right)
\end{align}
Hence the competitive equilibrium exists uniquely. 
\end{proof}}
{{The proof is provided in \addcite.}} 

We note that the competitive equilibrium in Theorem~\ref{mix_thrm_v2} does not need any extra requirements in market parameters, unlike the mixed market mechanism in subsection~\ref{subsec_2}. Additionally, the resulting market equilibrium with constant optimal participant bids can be solved fast enough for the needs of the real-time market using convex optimization. 

\section{Case Study}\label{sec_5}

In this section, we provide a numerical case study that analyses the mixed market mechanism based on decision-aware bidding from subsection~\ref{subsec_3}. For the day-ahead market, we formulate a two-day optimization horizon consisting of Day $1$ with $24$ time periods as the binding period and Day $2$ with $24$ periods as the advisory period. In the real-time market, we use a rolling time horizon window of $24$ h that has the 1st hour as the binding period and the rest as advisory periods. We assume that at each hour in the real-time market, the market realizes the updated demand only for the 1st (binding) interval. As the window rolls forward, the operator minimizes the dispatch cost subject to operational constraints for stage-wise and overall dispatch.  

We use forecast and real aggregate demand data for Aug 25-26, 2023, from the Millwood Zone in the New York ISO~\cite{nyisodata}. Furthermore, we assume one generator with aggregate cost coefficients $c = 20\$/(MW)^2$~\cite{matpower} and capacity limits $\underline{g} = 0,  \ \overline{g} = \max_t\{d_t\}$. Also, we assume one energy storage asset with variable capacity cost $B$ and variable capacity $E$. The empirical cost coefficient given by $bs := \rho B E$ where $\rho = 5.24 \times 10^{-4}$~\cite{shi2017optimal, bansal2020storage}. We assume a $4$-hour Lithium-ion battery such that the storage dispatch is bounded by $\underline{u} = -\frac{E}{4}$ and $\overline{u} = \frac{E}{4}$.

It is important to note that in real-time markets, the absence of a periodicity constraint may lead to periodicity constraint violations in cases of total two-stage dispatch. The operator only considers the optimal real-time dispatch for the first interval of the horizon as binding as the horizon window moves forward. Therefore, we use the underlying social planner problem~\eqref{social_planner} with and without periodicity constraints as a benchmark to compare the performance of the proposed mixed market mechanism from different perspectives.

\begin{figure}[ht]
    \centering
    {\includegraphics[width=\linewidth]{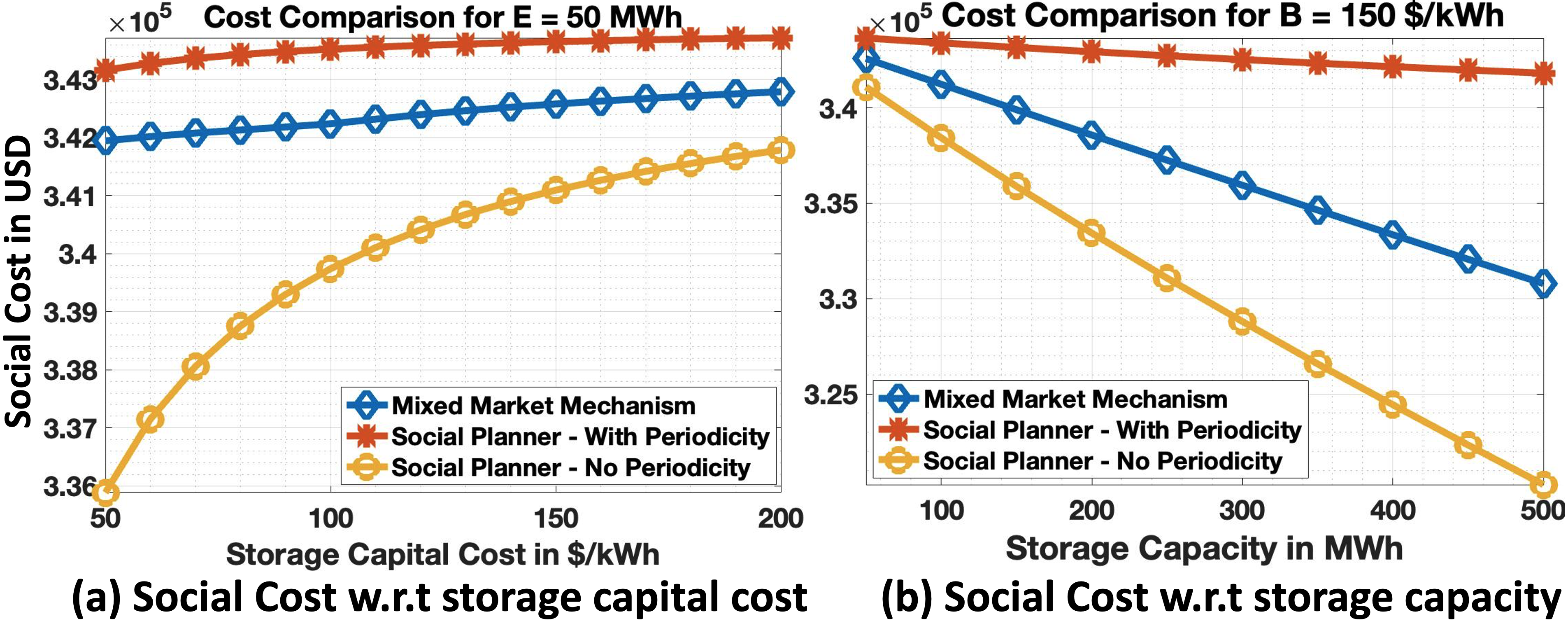} }
    \caption{Social cost in the proposed mixed market mechanism, social planner without periodicity constraints, and social planner with periodicity constraints w.r.t (left) storage capital cost and (right) storage capacity.}%
    \label{fig:cost}%
\end{figure}

We first consider the market perspective. Figure~\ref{fig:cost} illustrates the two-stage social cost for the three strategies, i.e., the proposed mixed market mechanism, the social planner with periodicity, and the social planner without periodicity, as we (a) change the storage capital cost for a fixed storage capacity of $E = 50 MWh$ and (b) change the storage capacity for a fixed storage capital cost $B = 150 \$/kWh$. The two-stage social cost of the proposed mechanism lies between the two bounds of the social planner problems in both the left and the right panels. This means that the proposed mechanism does not lead to significant losses in incentive compatibility. However, as we decrease the storage capital cost on the left panel or increase the storage capacity on the right panel, the difference in the social cost tends to increase. This suggests further study of the proposed mechanism in a more rigorous market setting.

\begin{figure}[ht]
    \centering
    {\includegraphics[width=\linewidth]{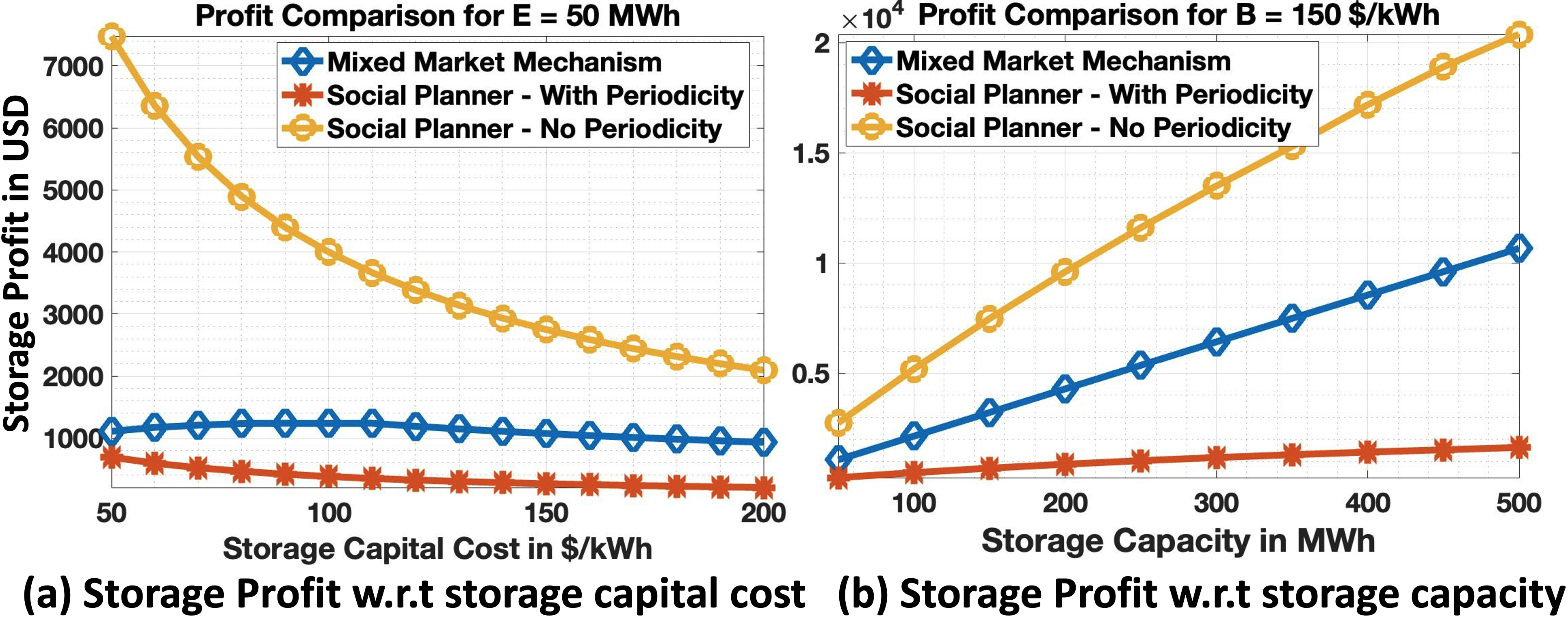} }
    \caption{Storage profit in the proposed mixed market mechanism, social planner without periodicity constraints, and social planner with periodicity constraints w.r.t (left) storage capital cost and (right) storage capacity.}%
    \label{fig:profit}%
\end{figure}

We next consider the resource owner perspective. Figure~\ref{fig:profit} illustrates the two-stage profit for the three strategies, i.e., the proposed mixed market mechanism, the social planner with periodicity, and the social planner without periodicity, as we (a) change the storage capital cost for a fixed storage capacity of $E = 50 MWh$ and (b) change the storage capacity for a fixed storage capital cost $B = 150 \$/kWh$. As expected, the energy storage profit also lies within the two bounds of the social planner problem. Similar to the previous case, the difference in profit between the mixed market mechanism and social planner without periodicity tends to increase as energy storage becomes cheaper or larger. However, note that profits tend to saturate at low capital costs, which results in greater utilization in real-time markets, as shown in the left panel.

\section{Concluding Remarks}\label{sec_6}

We study the market participation strategy of energy storage in two-stage markets while accounting for the temporally coupled storage operating cost in the form of its cycling cost. We propose two specific mixed market mechanisms where energy storage bids charge-discharge cycles in the day-ahead market followed by charge-discharge power bids in the real-time market. In the first mechanism, participants take into account their day-ahead decisions in their individual profit maximization problems. Alternatively, in the second approach, participants reflect on their day-ahead decisions in both bidding functions and individual problems. 

Although both market mechanisms result in a unique competitive equilibrium, the first requires an iterative best response algorithm which may not be desirable in real-time markets due to frequent market clearing. The second approach, on the other hand, can be implemented via convex programming. Numerical simulations for the real-world NYISO data show that the proposed mechanism lies within the two bounds of the social planner problem, implying a bounded loss in incentive compatibility. 

%\subsubsection*{ Future Extensions} It can be difficult to make quick decisions in the real-time market due to the complex and temporally coupled degradation cost function. However, our work is an important step towards designing a two-stage market for energy storage. While mixed market design shows promise, it still needs to be extended to include network constraints, opportunity costs at the end of the horizon, constraints on generator ramping, etc. We will leave these extensions to be addressed in future work.

\iffalse
Notes:
\begin{enumerate}
    \item DA - 2-day simulation; cost for each day, includes inequality and periodicity constraints
    \item Recall the optimal bid of storage in the real-time market is a function of total demand. As the window rolls forward, we have $D(1)$ as updated and the rest are from the forecast. We use this to find the optimal bid. ---- This is different from keeping a static D.
    \item No periodicity in real-time due to single binding interval 
\end{enumerate}

\fi

\bibliographystyle{IEEEtran}
\bibliography{cdc_24}

\begin{thebibliography}{10}
\providecommand{\url}[1]{#1}
\csname url@samestyle\endcsname
\providecommand{\newblock}{\relax}
\providecommand{\bibinfo}[2]{#2}
\providecommand{\BIBentrySTDinterwordspacing}{\spaceskip=0pt\relax}
\providecommand{\BIBentryALTinterwordstretchfactor}{4}
\providecommand{\BIBentryALTinterwordspacing}{\spaceskip=\fontdimen2\font plus
\BIBentryALTinterwordstretchfactor\fontdimen3\font minus \fontdimen4\font\relax}
\providecommand{\BIBforeignlanguage}[2]{{%
\expandafter\ifx\csname l@#1\endcsname\relax
\typeout{** WARNING: IEEEtran.bst: No hyphenation pattern has been}%
\typeout{** loaded for the language `#1'. Using the pattern for}%
\typeout{** the default language instead.}%
\else
\language=\csname l@#1\endcsname
\fi
#2}}
\providecommand{\BIBdecl}{\relax}
\BIBdecl

\bibitem{aeo2021}
``Annual energy outlook 2021 with projections to 2050.''\hskip 1em plus 0.5em minus 0.4em\relax Energy Information Administration, Washington, DC, 2021.

\bibitem{ferc2222}
``Participation of distributed energy resource aggregations in markets operated by regional transmission organizations and independent system operators, {Docket No.} rm18-9-000, {Order No.} 2222.''\hskip 1em plus 0.5em minus 0.4em\relax Department of Energy, {Federal} Energy Regulatory Commission, USA, 2020.

\bibitem{denholm2019potential}
P.~Denholm, J.~Nunemaker, P.~Gagnon, and W.~Cole, ``The {Potential} for battery energy storage to provide peaking capacity in the {United} {States},'' National Renewable Energy Laboratory (NREL), Tech. Rep., 2019.

\bibitem{solar_integr}
C.~A. Hill, M.~C. Such, D.~Chen, J.~Gonzalez, and W.~M. Grady, ``Battery energy storage for enabling integration of distributed solar power generation,'' \emph{IEEE Transactions on Smart Grid}, 2012.

\bibitem{wind_integ}
M.~A. Abdullah, K.~M. Muttaqi, D.~Sutanto, and A.~P. Agalgaonkar, ``An effective power dispatch control strategy to improve generation schedulability and supply reliability of a wind farm using a battery energy storage system,'' \emph{IEEE Transactions on Sustainable Energy}, vol.~6, no.~3, pp. 1093--1102, 2015.

\bibitem{trans_distr_supp}
R.~T. Elliott \emph{et~al.}, ``Sharing energy storage between transmission and distribution,'' \emph{IEEE Transactions on Power Systems}, 2019.

\bibitem{freq_reserve_req}
D.~Pudjianto, M.~Aunedi, P.~Djapic, and G.~Strbac, ``Whole-systems assessment of the value of energy storage in low-carbon electricity systems,'' \emph{IEEE Transactions on Smart Grid}, vol.~5, no.~2, pp. 1098--1109, 2014.

\bibitem{bansal_epri}
\BIBentryALTinterwordspacing
N.~G. Singhal, R.~K. Bansal, J.~M. Kemp, E.~Ela, and M.~Heleno, ``Integration of hybrids into wholesale power markets,'' Lawrence Berkeley National Laboratory (LBNL), Berkeley, CA (United States), Tech. Rep., 8 2023. [Online]. Available: \url{https://www.osti.gov/biblio/1994813}
\BIBentrySTDinterwordspacing

\bibitem{he2015optimal}
G.~He, Q.~Chen, C.~Kang, P.~Pinson, and Q.~Xia, ``Optimal bidding strategy of battery storage in power markets considering performance-based regulation and battery cycle life,'' \emph{IEEE Transactions on Smart Grid}, vol.~7, no.~5, pp. 2359--2367, 2015.

\bibitem{dheepak_model_1}
D.~Krishnamurthy, C.~Uckun, Z.~Zhou, P.~R. Thimmapuram, and A.~Botterud, ``Energy storage arbitrage under day-ahead and real-time price uncertainty,'' \emph{IEEE Transactions on Power Systems}, 2018.

\bibitem{ramteen_model_1}
S.~Bhattacharjee, R.~Sioshansi, and H.~Zareipour, ``Energy storage participation in wholesale markets: The impact of state-of-energy management,'' \emph{IEEE Open Access Journal of Power and Energy}, vol.~9, pp. 173--182, 2022.

\bibitem{bansal2021market}
R.~K. Bansal, P.~You, D.~F. Gayme, and E.~Mallada, ``A market mechanism for truthful bidding with energy storage,'' \emph{Electric Power Systems Research}, vol. 211, 2022.

\bibitem{BXU_model_2}
N.~Zheng, X.~Qin, D.~Wu, G.~Murtaugh, and B.~Xu, ``Energy storage state-of-charge market model,'' \emph{IEEE Transactions on Energy Markets, Policy and Regulation}, vol.~1, no.~1, pp. 11--22, 2023.

\bibitem{Cong_model_2}
C.~Chen, S.~Li, and L.~Tong, ``Multi-interval energy-reserve co-optimization with soc-dependent bids from battery storage,'' 2024.

\bibitem{bansal2020storage}
R.~K. Bansal, P.~You, D.~F. Gayme, and E.~Mallada, ``Storage degradation aware economic dispatch,'' in \emph{2021 American Control Conference (ACC)}, 2021, pp. 589--595.

\bibitem{nyisodata}
``{ISO-New York} : Energy market and operational data,'' Available at \url{https://www.nyiso.com/custom-reports/} (2023).

\bibitem{matpower}
R.~D. Zimmerman and C.~E. {Murillo-Sanchez}, ``Matpower (version 7.0) [software].'' Available at \url{https://matpower.org/} (2019).

\bibitem{shi2017optimal}
Y.~Shi, B.~Xu, Y.~Tan, D.~Kirschen, and B.~Zhang, ``Optimal battery control under cycle aging mechanisms in pay for performance settings,'' \emph{IEEE Transactions on Automatic Control}, vol.~64, no.~6, pp. 2324--2339, 2019.

\end{thebibliography}

\end{document}